\documentclass[11pt]{article}
\usepackage{graphicx}
\usepackage{amssymb,amsfonts,amsthm}

\usepackage{amsmath,amsthm,amssymb}
\usepackage{amscd}
\usepackage{amsfonts}
\usepackage{color}
\newcommand{\la}{\langle}
\newcommand{\ra}{\rangle}

\newtheorem{theorem}{Theorem}
\newtheorem{lemma}[theorem]{Lemma}

\theoremstyle{definition}

\newtheorem{remark}[theorem]{Remark}

\newcommand {\iv}{^{-1}}

\begin{document}

\title{A 2-generated 2-related group with no non-trivial finite factors}
 \author{A.Yu. Ol'shanskii, M.V. Sapir\thanks{The
authors were supported in part by the NSF grants DMS 0245600 and DMS
0455881.  In addition, the research of the first author was
supported in part by the Russian Fund for Basic Research grant
05-01-00895, and the research of the second author was supported by
a BSF (USA-Israeli) grant.}}
\date{}
\maketitle


\begin{abstract}
We construct a 2-generated 2-related group without non-trivial
finite factors. That answers a question of J. Button.
\end{abstract}

Problem 1.12 from \cite{KT} (attributed to  Magnus and included in
the 1965 edition of \cite{KT} by Greendlinger) asks whether the
triviality problem for groups given by balanced presentations with
$n\ge 2$ generators is decidable. That problem would have an easy
solution if every non-trivial $n$-generated $n$-related group would
have a non-trivial finite factor. Indeed, if that was the case, then
in order to check triviality of that group one could simultaneously
list its finite factors and all the relations of the group. The
group would be non-trivial if it has a non-trivial finite factor and
trivial if one can deduce relations $x=1$ for all generators $x$
(see more details in \cite[Section 2.6]{KS}). Unfortunately, it is
well known that there are infinite groups with balanced
presentations which do not have non-trivial finite factors. In
particular, such is Higman's group $\la a,b,c,d\mid a^b=a^2,
b^c=b^2, c^d=c^2, d^a=d^2\ra$ \cite{Higman}. But for groups having
balanced presentations with fewer than 4 generators the answer seems
to be not present in the literature. The question whether a
2-generator 2-relator group without non-trivial finite factors
exists was asked by Jack Button (we found out about this question
from Ian Leary). In this note, we give a positive answer to this
question. Our example also has a 3-generated balanced presentation,
of course (add a new generator $x$ and a relation $x=1$).

\begin{theorem}\label{m} There exists an infinite
2-generated 2-related group with no non-trivial finite factors.
\end{theorem}

\proof Let $G$ be the Baumslag group $\la a,t\mid a^{a^t}=a^2\ra$
where $x^y$ denotes $y\iv xy$. It is known \cite{Ba} that all finite
factors of $G$ are cyclic. Moreover, the image of $a$ is $1$ in
every finite factor of $G$.

It will be convenient for us to represent $G$ as an HNN extension
$\la a,b,t\mid a^t=b, a^b=a^2\ra$ of the Baumslag-Solitar group
$H$=$\la a,b\mid a^b=a^2\ra$ with associated subgroups $A=\la
a\ra$ and $B=\la b\ra$ such that $t^{-1}At=B$ in $G$.

Recall \cite{LS} if an element $g\in G$ is equal to the product
$h_0t^{\epsilon_1}h_1\dots t^{\epsilon_n}h_n$  where $h_0,\dots,
h_n\in H$ and $\epsilon_i=\pm 1$ ($i=1,\dots,n$), then the
$t$-length of $g$ is $n$. The product is called {\it reduced} if it
has  no
 occurrences of $t^{-1}xt$ with $x\in A$ or $tyt^{-1}$ with $y\in B$.
 It is called {\em cyclically reduced} if every cyclic permutation
 of that product is reduced.
 By Britton's lemma \cite{LS} two reduced products representing the same element in $G$
 have equal $t$-lengths. The following property holds for the Baumslag group $G$.

\medskip

($\star$) {\it For every integer $n\ge 2$ and $x,y\in H$, we have
$t^n yt^{-n}=x$ in $G$ if and only if $x=y=1$.}

\medskip

Indeed, the $t$-length of the right-hand side is $0$, and so $y\in
B$, and $tyt^{-1}=z$ in $G$, where $z\in A$. But also $z\in B$
since $tzt^{-1}=x$ in $G$ . Hence $z\in A\cap B=\{1\}$, and
$x=y=1$.

 Let $r$ be any word of the form

$$bt^{u_1}a t^{-u_2}b t^{u_3}...at^{-u_l}$$ where the numbers $l$ and $u_i$ satisfy the following conditions:

\begin{itemize}

\item[(1)] $l$ is even, and $u_i$ are different integers, $u_i\ge
2$

 \item[(2)] the
total exponent of $t$, i.e. $u_1-u_2+u_3...-u_l$ is $1$

\item [(3)] $\max (u_1+u_2+u_3, u_2+u_3+u_4,\dots,u_{l-1}+u_l+u_1,
u_l+u_1+u_2)< \frac16 \sum_{i=1}^l u_i$

\end{itemize}
(for example one can take $l=20$, $u_{2i-1}=100+2i-1$ for
$i=1,\dots,9$, $u_{19}=130$, $u_{2i}=100+2i$ for $i=1,\dots,10$).
Consider the factor-group $K=\la G\mid r=1\ra$. The images of $a$
and $b=a^t$ vanish in every finite factor-group $F$ of $K$ since the
same property holds for $G$. It follows from (2) and equality $r=1$
that the image of $t$ in $F$ is also trivial. Therefore $K$ does not
have non-trivial finite factors, and it remains to prove that $K$
itself is not trivial.

The product $r$ is cyclically reduced. Let $R$ be the set of all
cyclically reduced forms of the conjugate elements of $r$ and
$r^{-1}$ in $G$.

\begin{lemma} \label{Coll} (D.Collins \cite{LS}
, IV.2.5)

Let $w=h_1t^{\epsilon_1}\dots h_nt^{\epsilon_n}$ ($n\ge 1$)
 and $w'$ be conjugate
cyclically reduced elements in an HNN extension $H^{\star}$ of a
group $H$ ( $h_1,\dots, h_n\in H$). Then $w'$ is equal in
$H^{\star}$ to $h^{-1}w^{\star}h$ for a cyclic permutation
$w^{\star}$ of $w$ and some $h\in H$.

\end{lemma}

\begin{lemma} \label {equal} Let $r_1$ and $r_2$ belong to $R$ and
equal in $G$ to some reduced products starting with $t^{\pm
u_s}xt^{\mp u_{s+1}}x'$ and $t^{\pm u_s}yt^{\mp u_{s+1}}y'$,
respectively, where $x, x', y$, and $y'$ are nontrivial elements
of $H$, and the subscripts are taken modulo $l$. Then $r_1=r_2$ in
$G$.
\end{lemma}

\proof By Lemma \ref{Coll} and condition (1), $r_1$ and $r_2$ are
conjugates in $G$ of the same cyclic permutation of $r$ by some
elements of $H$. Hence $r_2=hr_1h^{-1}$ in $G$ for some $h\in H$.
But the factors $t^{\pm u_s}$ must cancel out in the product of
$r_1^{-1}r_2$ by the lemma condition. Therefore they also must
cancel in the product of reduced $r_1$ and $hr_1h^{-1}$. Hence
$\dots t^{\mp u_s}ht^{\pm u_s}\in H$. It follows from Condition
(1) and Property ($\star$) that $h=1$, and therefore $r_1=r_2$.
\endproof

\begin{lemma} \label{main} Let $r_1=vw_1$ and $r_2=vw_2$ in $G$,
where $r_1,r_2\in R$, and these products are reduced in $G$.
Assume that a reduced form of $v$ starts with $t^uh_1t^{\pm
u_s}h_2t^{\mp u_{s+1}}h_3$ for some non-trivial $h_1,h_2,h_3\in H$
and an integer $u$ (subscripts are taken modulo $l$). Then
$r_1=r_2$ in $G$.
\end{lemma}

\proof Denote the prefix $t^uh_1$ of $v$ by $z$. There are reduced
forms of $z^{-1}r_1z$ and of $z^{-1}r_2z$ both started with
$t^{\pm u_s}h_2t^{\mp u_{s+1}}h_3$. By Lemma \ref{equal}, we have
$z^{-1}r_1z=z^{-1}r_2z$, and so $r_1=r_2$ in $G$.\endproof

\medskip

{\bf Proof of the theorem.} Let $r_1\in R$, and assume that $r_1$ is
conjugate of $r$ in $G$. (The proof is similar, if it is conjugate
of $r^{-1}$.) If $r_1$ and a word $r_2\in R$ have a left piece $v$
as in the formulation of Lemma \ref{main}, then $r_1=r_2$ by this
lemma. Otherwise $v$ contains at most $3$ $t$-blocks, and so its
$t$-length is less than $\frac16|r_1|$  by Condition (3). Hence our
presentation $K=\la G\mid r=1\ra$ satisfies the small cancellation
condition $C'(\frac16)$ by definition (see \cite{LS}, V.11). Then by
Theorem 11.6 \cite{LS}, the canonical homomorphism of $H$ into $K$
is injective. Therefore the group $K$ is infinite, as required.

\begin{remark} Note that not only finite, but all torsion homomorphic images
of the group $K$ are trivial because the image of $a$ is trivial in
every torsion factor-group of $G$ \cite{Ba}. Since every hyperbolic
group is residually torsion (\cite{Gr},\cite{Ol}), it follows that
the image of $K$ under any homomorphism of $K$ into a hyperbolic
group is trivial. On the other hand, $K$ has continuously many
normal subgroups, and moreover, it is $SQ$-universal. To prove this,
one can just combine the above small cancellation argument with the
constructions from the proofs of theorems 11.3 and 11.7 in
\cite{LS}.

\end{remark}

\noindent Alexander Yu. Olshanskii:\\
{\small \sc Department of Mathematics, Vanderbilt University , Nashville, TN 37240.\\
Department of Mathematics, Moscow State University, Moscow, 119899,
Russia.\\} {\it E-mail:} {\tt alexander.olshanskiy@vanderbilt.edu}

\vspace{3mm}

\noindent Mark V. Sapir: {\small\sc \\ Department of Mathematics,
Vanderbilt University, Nashville, TN 37240.\\} {\it E-mail: } {\tt
m.sapir@vanderbilt.edu}


\begin{thebibliography}{xxxxx}

\addcontentsline{toc}{section}{Bibliography}

\bibitem[Ba]{Ba} G. Baumslag. A non-cyclic one-relator group all of whose
finite quotients are cyclic.  J. Austral. Math. Soc.  10 (1969)
497--498.

\bibitem[Gr]{Gr} M.Gromov. Hyperbolic groups, in Essays in Group
Theory (S.M.Gersten, ed.), M.S.R.I. Pub. 8, Springer (1987),
75-263.

\bibitem[Hi]{Higman} G. Higman. A finitely generated infinite simple group, J. London
Math. Soc., 26, (1951) 61--64.

\bibitem[KS]{KS} O. Kharlapovich, M. Sapir.
Algorithmic problems in varieties, Internat. J. Algebra Comput. 5
(1995), no. 4-5, 379--602.

\bibitem[KT]{KT} Kourov notebook. Unsolved problems in the theory of
groups.
Academy of Sciences of the USSR. Siberian Branch. Institute of
Mathematics Izdat. Sibirsk. Otdel. Akad. Nauk SSSR, Novosibirsk 1965
18 pp.


\bibitem[LS]{LS}
R.C. Lyndon, P.E. Shupp, Combinatorial Group Theory,
Springer--Verlag, 1977.



\bibitem[Ol]{Ol}
A.Yu.Olshanskii. On residualing homomorphisms and G-subgroups of
hyperbolic groups, Internat. J. of Algebra and Comput., 3 (1993),
1-44.



\end{thebibliography}
\end{document}